\documentclass{article}
\usepackage{spconf,amsmath,graphicx,hyperref}
\usepackage{epstopdf}
\usepackage{amsmath}
\usepackage{amssymb}
\usepackage{amsthm}
\usepackage{color}
\usepackage{enumitem}
\usepackage{hyperref}
 \hypersetup{
     colorlinks=true,
     linkcolor=black,
     filecolor=black,
     citecolor =black,      
     urlcolor=black,
     }
\usepackage{bm}
\usepackage{mathrsfs}
\usepackage[scr=boondox]{mathalpha}
 \usepackage{booktabs} 
\usepackage{extarrows}
\usepackage{graphicx}
\usepackage{subfig}
\usepackage{algorithmic,algorithm}
\usepackage{mathrsfs}
\usepackage{cmll}
\usepackage{mathtools}
\newtheorem{theorem}{Theorem}
\newtheorem{proposition}{Proposition}
\newtheorem{lemma}{Lemma}

\theoremstyle{definition}

\newtheorem{definition}{Definition}
\newtheorem{assumption}{Assumption}

\theoremstyle{remark}

\newcommand{\BE}{\begin{equation}}
\newcommand{\EE}{\end{equation}}
\newcommand{\BS}{\begin{subequations}}
\newcommand{\ES}{\end{subequations}}
\newcommand{\UT}{\mathsf{T}}   %% upright T

\newcommand{\wc}{\xrightarrow{\mathcal{W}_2}} % Wasserstein-2 convergence
\newcommand{\ac}{\xrightarrow{a.s.}} % Almost sure convergence

\providecommand{\bm}{\boldsymbol{}}

\title{Orthogonal Approximate Message Passing Algorithms for Rectangular Spiked Matrix Models with Rotationally Invariant Noise}
\name{Haohua Chen$^{*}$ \qquad Songbin Liu$^{\dagger}$ \qquad Junjie Ma$^{*}$}

\address{$^{*}$Academy of Mathematics and Systems Science, Chinese Academy of Sciences, Beijing, China \\
         $^{\dagger}$Department of Electrical Engineering, Columbia University, New York, USA \\
         Email: chenhaohua25@mails.ucas.ac.cn, sl5878@columbia.edu, majunjie@lsec.cc.ac.cn
}
\begin{document}
\ninept
\maketitle
\begin{abstract}
We propose an orthogonal approximate message passing (OAMP) algorithm for signal estimation in the rectangular spiked matrix model with general rotationally invariant (RI) noise. We establish a rigorous state evolution that exactly characterizes the high-dimensional dynamics of the algorithm. Building on this framework, we derive an optimal variant of OAMP that minimizes the predicted mean-squared error at each iteration. For the special case of i.i.d. Gaussian noise, the fixed point of the proposed OAMP algorithm coincides with that of the standard AMP algorithm. For general RI noise models, we conjecture that the optimal OAMP algorithm is statistically optimal within a broad class of iterative methods, and achieves Bayes-optimal performance in certain regimes.
\end{abstract}
\section{Introduction}\label{sec:intro}
We consider the estimation of rank-one signals \(\bm{u}_*\in \mathbb{R}^M\), \(\bm{v}_* \in \mathbb{R}^N\) from a \emph{rectangular} spiked matrix model:
\begin{equation}\label{eq:rectangular_spiked_model}
\bm{Y} = \frac{\theta}{\sqrt{MN}} \bm{u}_* \bm{v}_*^\top + \bm{W} \in \mathbb{R}^{M\times N},
\end{equation}
where $\bm{W} \in \mathbb{R}^{M \times N}$ is a rotationally invariant (RI) noise matrix and \(\theta\geq 0\) is a signal-to-noise ratio (SNR) parameter. The spiked model is a powerful tool for analyzing data where the number of features is of a comparable scale to the number of samples, with notable applications in financial data analysis \cite{Bouchaud2009Financial}, community detection \cite{Abbe2016block}.

For the spiked model in \eqref{eq:rectangular_spiked_model} with i.i.d. Gaussian noise, a well-established phase transition dictates a sharp signal-to-noise ratio (SNR) threshold for signal recovery via Principal Component Analysis (PCA) \cite{Feral2007Deformation, Bai2010Spectral}. It is known that the principal eigenvector aligns non-trivially with the underlying signal if and only if the SNR exceeds this critical threshold. The performance of standard PCA can be further improved by incorporating prior structural information about the signal (e.g., sparsity) \cite{Zou2006Sparse} --\nocite{Deshpande2014PCA,Lesieur2015PCA,Montanari2016nonnegative} \cite{Lesieur2017lowrank}. In this context, the approximate message passing (AMP) algorithm, originally developed for compressed sensing  \cite{Donoho2009cs} --\nocite{Bayati2011dense,Ma2017OAMP,Rangan2019VAMP} \cite{Takeuchi2020Rigor}, offers the provably optimal performance among generalized first order methods \cite{montanari2021estimation}, and attains the minimum mean-squared error (MMSE) predicted by statistical physics in certain regimes \cite{Barbier2016SymmetricReplica} -- \nocite{Miolane2018Non-symmetricReplica,Alaoui2018SpikedReplica}\cite{Lelarge2019Fundamental}.

For the rectangular spiked model with rotationally invariant (RI) noise, the PCA phase transition is also well understood \cite{Florent2012singular}, but its performance is generally suboptimal when priors are available. To address this, AMP-type algorithms for RI models (referred to as {RI-AMP} here) have been proposed \cite{fan2022ROT,zhong2024spectral}, though their MSE performance optimality are not established \cite{fan2022ROT}. For the \textit{symmetric} counterpart of \eqref{eq:rectangular_spiked_model}, a major advance was achieved in \cite{Dudeja2024optimality}. Building on the OAMP framework \cite{Ma2017OAMP} and inspired by insights from \cite{barbier2023fundamental}, \cite{Dudeja2024optimality} demonstrated that optimality within a broad class of algorithms requires not only prior-informed iterative denoising but also a pre-processing of the observation matrix $\bm{Y}$ based on the spectrum of the noise matrix. This pre-processing connects naturally to classical results in the literature of covariance matrix shrinkage and matrix denoising \cite{ledoit2012nonlinear} -- \nocite{Nada2014Optshrink,Bun2016SymmetricRIE}\cite{Pourkamali2025RectRIE}. Furthermore, the performance of the OAMP algorithm in \cite{Dudeja2024optimality} was shown to be optimal among a large class of algorithms and match the replica predictions \cite{barbier2023fundamental,Barbier2025TAP} under certain conditions.

In this work, we extend the OAMP framework \cite{Dudeja2024optimality} to the rectangular spiked matrix model. We develop a tractable state evolution (SE) for the proposed OAMP algorithm and derive the corresponding optimal denoisers. We further show that the fixed-point equations of the OAMP algorithm recover known results in the i.i.d. Gaussian case \cite{montanari2021estimation}, thereby supporting the conjecture of its potential optimality. 

Detailed analyses and full proofs of our main results are provided in the extended version of this work \cite{Chen2025rect}.

\subsection{Notation and Assumptions.} 
Let $(\bm{x}_1, \dots, \bm{x}_K)$ be a collection of vectors, where each vector $\bm{x}_k = (x_{k,1}, \dots, x_{k,d})^\top \in \mathbb{R}^d$. The convergence of this collection to a tuple of random variables $(\mathsf{X}_1, \dots, \mathsf{X}_K)$ in the Wasserstein-2 distance, denoted $(\bm{x}_1, \dots, \bm{x}_K) \wc (\mathsf{X}_1, \dots, \mathsf{X}_K)$, is defined by the condition that for all pseudo-Lipschitz test functions $h: \mathbb{R}^K \to \mathbb{R}$ of order 2, the limit $\lim_{d \to \infty} \frac{1}{d} \sum_{i=1}^d h(x_{1,i}, x_{2,i}, \dots, x_{K,i}) \xrightarrow{a.s.} \mathbb{E}[h(\mathsf{X}_1, \mathsf{X}_2, \dots, \mathsf{X}_K)]$ holds almost surely. Inner product with respect to any finite measure $\chi(\lambda)$ on $\mathbb{R}$ for any bounded, Borel measurable functions $p: \mathbb{R} \mapsto \mathbb{R}$ is denoted by $\langle p(\lambda)  \rangle_{\chi} \triangleq \int p(\lambda) d\chi(\lambda)$. The Stieltjes transform of a finite measure $\mu$ is defined as
\begin{equation*}
\mathcal{S}_\mu(z) \triangleq \int \frac{1}{z - \lambda}\mathrm{d}\mu(\lambda),\quad\forall z \in \mathbb{C}\setminus\text{supp}(\mu).
\end{equation*}

We impose the following assumptions throughout this paper.

\begin{assumption}[Signal and Noise Models]\label{assump:main}
$\ $
\begin{enumerate}
\item[(a)] We consider the high-dimensional limit where $M, N \to \infty$ with an aspect ratio $M/N \to \delta \in (0,1]$. 
    \item[(b)] The signal and their side information converge as $(\bm{u}_*, \bm{a}) \wc \pi_U$ and $(\bm{v}_*, \bm{b}) \wc \pi_V$, and are normalized to unit variance such that $\|\bm{u}_*\|^2/M \xrightarrow{a.s.} 1$ and $\|\bm{v}_*\|^2/N \xrightarrow{a.s.} 1$.

    \item[(c)] The noise matrix is given by $\bm{W} = \bm{U}\mathrm{diag}(\bm{\sigma})\bm{V}^\UT \in \mathbb{R}^{M\times N}$, where $\bm{U}$ and $\bm{V}$ are independent Haar-distributed orthogonal matrices. We assume $\|\bm{WW}^\UT\|_{\mathrm{op}} \leq C$ for some dimension-independent constant $C$. The empirical spectral distribution of the matrix $\bm{W}\bm{W}^\UT$ converges weakly to a deterministic measure $\mu$, which is absolutely continuous with a Hölder continuous density on a compact support. The limiting spectral measure of $\bm{W}^\UT \bm{W}$ is then $\tilde{\mu} \triangleq \delta\mu+(1-\delta)\delta_{\{0\}}$.
\end{enumerate}
\end{assumption}

\section{Preliminaries on Spectral Analysis}\label{sec: Spec}

We begin with some spectral analysis results that will be used in the main results (Section \ref{sec:algo}) of this paper.

\begin{definition}[Signal–Eigenspace Spectral Measures]\label{def:spectral_measures} 
Let $(\lambda_i(\cdot), \bm{u}_i(\cdot))$ denote the eigenvalue–eigenvector pairs of a matrix.
\begin{enumerate}
    \item[(a)] Define the following weighted empirical spectral distributions:
    \begin{align*}
        \nu_{M,1} &\triangleq \frac{1}{M}\sum_{i=1}^M 
        \langle \bm{u}_i(\bm{Y}\bm{Y}^\top), \bm{u}_* \rangle^2 
        \, \delta_{\lambda_i(\bm{Y}\bm{Y}^\top)}, \\
        \nu_{N,2} &\triangleq \frac{1}{N}\sum_{i=1}^N 
        \langle \bm{u}_i(\bm{Y}^\top \bm{Y}), \bm{v}_* \rangle^2 
        \, \delta_{\lambda_i(\bm{Y}^\top \bm{Y})}.
    \end{align*}

    \item[(b)] Consider the symmetric dilation of $\bm{Y}$:
    \[
        \widehat{\bm{Y}} \triangleq
        \begin{bmatrix} 
            \bm{0} & \bm{Y} \\ 
            \bm{Y}^\top & \bm{0}
        \end{bmatrix}
        \in \mathbb{R}^{L \times L}, 
        \qquad L \triangleq M+N.
    \]
    Define $\widehat{\bm{u}}_* \triangleq [\bm{u}_*^\top, \bm{0}^\top]^\top$ 
    and $\widehat{\bm{v}}_* \triangleq [\bm{0}^\top, \bm{v}_*^\top]^\top$.  
    Define an associated empirical signed spectral measure:
    \[
        \nu_{L,3} \triangleq \frac{1}{L}\sum_{i=1}^L 
        \langle \bm{u}_i(\widehat{\bm{Y}}), \widehat{\bm{u}}_* \rangle 
        \langle \bm{u}_i(\widehat{\bm{Y}}), \widehat{\bm{v}}_* \rangle 
        \, \delta_{\lambda_i(\widehat{\bm{Y}})}.
    \]
\end{enumerate}
\end{definition}

Three \textit{shrinkage functions} $\varphi_1$, $\varphi_2$ and $\varphi_3$, which are related to the limit of the above empirical measures, will be important for our OAMP algorithm. To state the result, define:
    \begin{equation}\label{eq:Ctrans}
\mathcal{C}(z) \triangleq z \mathcal{S}_\mu(z) \left(\delta \mathcal{S}_\mu(z) + (1-\delta)/z\right),\quad \forall z\in\mathbb{C}\backslash\text{supp}(\mu),
\end{equation}
where $\mathcal{S}_\mu$ is the Stieltjes transform of the limiting spectrum of the noise matrix $\mu$; see Assumption \ref{assump:main}.  The functions $\varphi_i$ are defined as:
\begin{subequations}\label{Eqn:varphi_functions}
    \begin{align}
    \varphi_1(\lambda) &\triangleq \frac{1 + \delta\theta^2\pi^2\lambda\left(\mathcal{H}_{\mu}(\lambda)^2 + \mu(\lambda)^2\right)}{\lim_{\epsilon \to 0^+} \left|1-\theta^2\mathcal{C}(\lambda-i\epsilon)\right|^2}, \label{eq:phiu} \\
    \varphi_3(\lambda) &\triangleq \frac{\theta\left(1-\delta + 2\delta\pi\lambda\mathcal{H}_{\mu}(\lambda)\right)}{\lim_{\epsilon \to 0^+} \left|1-\theta^2\mathcal{C}(\lambda-i\epsilon)\right|^2} \cdot \mathbf{1}_{\{\lambda \neq 0\}}, \label{eq:phiuv} \\
    \varphi_2(\lambda) &\triangleq 
    \begin{cases}
        \delta \varphi_1(\lambda) + \frac{\theta(1-\delta)}{\lambda}\varphi_3(\lambda), & \text{if } \lambda > 0, \\
        \frac{\delta}{1-\theta^2 (1-\delta)\pi\mathcal{H}(0)}, & \text{if } \lambda = 0,
    \end{cases}
    \label{eq:phiv}
\end{align}
\end{subequations}
where, $\mu(\lambda)$ denotes the density function of the measure $\mu$ and $\mathcal{H}_{\mu}$ denotes its Hilbert transform:\begin{equation}
\mathcal{H}_{\mu}(x) \triangleq \frac{1}{\pi} \text{P.V.} \int \frac{1}{x - \lambda} \mu(\lambda)\, \mathrm{d}\lambda.
\end{equation}
The functions $\varphi_i$ can be expressed explicitly using $\mu(\lambda)$ and its Hilbert transform by applying the Sokhotski–Plemelj formula~\cite{Blanchard2015Hilbert}:
\begin{equation}\label{eq:denominator_expanded}
\begin{split}
    &\lim_{\epsilon \to 0^+} |1-\theta^2\mathcal{C}(\lambda-i\epsilon)|^2 \\
    &\quad= \Big\{1 - \theta^2\left[ (1-\delta)\pi\mathcal{H}_{\mu}(\lambda) - \delta\pi^2\lambda ( \mu^2(\lambda) - \mathcal{H}_{\mu}^2(\lambda) ) \right] \Big\}^2 \\
    &\qquad + \Big\{\pi\theta^2\mu(\lambda)\left[1-\delta + 2\delta\pi\lambda\mathcal{H}_{\mu}(\lambda)\right] \Big\}^2.
\end{split}
\end{equation}

The following lemma characterizes the weak limit of the spectral measures in Definition~\ref{def:spectral_measures}. An analogous statement for the symmetric setting appears in \cite[Lemma 1]{Dudeja2024optimality}.

\begin{lemma}\label{lem:spectral_measures}
The following hold:
\begin{enumerate}  
\item [(a)] The measures $\nu_{M,1}$, $\nu_{N,2}$, $\nu_{L,3}$ in Definition \ref{def:spectral_measures} converge weakly almost surely to deterministic, compactly supported measures $\nu_1$, $\nu_2$ and $\nu_3$. Further, $\nu_1$ and $\nu_2$
are probability measures on $\mathbb{R}_+$, while $\nu_3$ is a signed measure on $\mathbb{R}$. 
\item  [(b)] For $z \in \mathbb{C}\setminus\mathbb{R}$, the Stieltjes transforms of $\nu_1$, $\nu_2$ and $\nu_3$ are given by:
    \begin{align}
        &\mathcal{S}_{\nu_1}(z) = \frac{\mathcal{S}_{\mu}(z)}{1-\theta^2 \mathcal{C}(z)}, \quad
        \mathcal{S}_{\nu_2}(z) = \frac{\delta \mathcal{S}_{\mu}(z) +\frac{1-\delta}{z} }{1-\theta^2 \mathcal{C}(z)}, \\& \mathcal{S}_{\nu_3}(z) = \frac{\sqrt{\delta}}{1+\delta}\cdot\frac{\theta \mathcal{C}(z^2)}{1- \theta^2 \mathcal{C}(z^2)}.
    \end{align}
    \item[(c)] Let $\nu_i=\nu_i^\parallel + \nu_i^\perp$ be the Lebesgue decomposition of $\nu_i$ into an absolutely continuous part $\nu_i^\parallel$ and a singular part $\nu_i^\perp$, for $i\in\{1,2,3\}$. The densities of the absolutely continuous parts are given by:
        \begin{align}
        \frac{d\nu_{i}^{\parallel}}{d\lambda} &= \mu(\lambda)\varphi_i(\lambda), \quad \text{for } i \in \{1,2\}, \label{eq:density_nu12}\\
        \frac{d\nu^{\parallel}_3}{d\sigma} &=  \frac{\sqrt{\delta}}{1+\delta}\mathrm{sign}(\sigma)\cdot\mu(\sigma^2)\varphi_{3}(\sigma^2).\label{eq:density_nu3}
        \end{align}
 
    \item [(d)] The measures have point masses at the roots of 
   \begin{equation}
   1-\theta^2 \mathcal{C}(\lambda)=0,\quad \text{for }\lambda \in\mathbb{R}\backslash \mathrm{supp}(\mu).
   \end{equation}
Let $\lambda_*>0$ be such a non-negative root. The corresponding masses are:
        \begin{align*}
         \nu_1(\{\lambda_*\})&= \frac{\mathcal{S}_{\mu} (\lambda_*)}{-\theta^2\mathcal{C}'(\lambda_*)},  \quad \nu_2(\{\lambda_*\})= \frac{\delta \mathcal{S}_{\mu} (\lambda_*) + \frac{1-\delta}{\lambda_*}}{-\theta^2\mathcal{C}'(\lambda_*)}, \\
         \nu_3(\{\pm\sigma_*\})&= \mp\frac{\sqrt{\delta}}{1+\delta}\frac{1}{2\theta^3\sigma_*\mathcal{C}'(\sigma^2_*)}, \quad \forall \sigma_* = \sqrt{\lambda_*} >0.
        \end{align*}
        Additionally, for $\delta<1$, $\nu_2$ exhibits a point mass at the origin:
        \begin{align}
            \nu_2(\{0\}) &= \frac{1-\delta}{1-\theta^2 (1-\delta)\pi\mathcal{H}(0)}. \label{eq:point_mass_zero}
        \end{align}
\end{enumerate}
\end{lemma}

\section{Orthogonal AMP Algorithm}\label{sec:algo}

The main results are presented in this section.

\subsection{OAMP Algorithm}
We propose a class of orthogonal approximate message passing (OAMP) algorithms for the rectangular model \eqref{eq:rectangular_spiked_model}. This class extends the OAMP algorithm introduced for symmetric spiked matrix models in \cite{Dudeja2024optimality}.

\begin{definition}[OAMP Algorithms]\label{def:OAMP}
An OAMP algorithm generates sequences of iterates $\{\bm{u}_t\}$ and $\{\bm{v}_t\}$ via the update rules:
\begin{align}
 &\bm{u}_t = F_{t}(\bm{Y}\bm{Y}^\UT) f_{t}(\bm{u}_{\leq t-1}; \bm{a}) +\widetilde{F}_{t}(\bm{Y}\bm{Y}^\UT)\bm{Y} g_{t}(\bm{v}_{\leq t-1}; \bm{b}) , \label{eq:OAMP Algo u}\\
 &\bm{v}_t = G_{t}(\bm{Y}^\UT\bm{Y}) g_{t}(\bm{v}_{\leq t-1}; \bm{b})+\widetilde{G}_{t}(\bm{Y}^\UT\bm{Y})\bm{Y}^\UT f_{t}(\bm{u}_{\leq t-1}; \bm{a}), \notag
\end{align}
where $f_t, g_t$ are entry-wise \textit{iterate denoisers} that depend on the history of past iterates $\bm{u}_{\leq t-1} \triangleq (\bm{u}_1, \dots, \bm{u}_{t-1})$, $\bm{v}_{\leq t-1} \triangleq (\bm{v}_1, \dots, \bm{v}_{t-1})$ and the side information $\bm{a}, \bm{b}$. The \textit{matrix denoisers} $F_t, \tilde{F}_t, G_t, \tilde{G}_t$ are scalar functions applied to the eigenvalues without changing the eigenvectors. For example, let the eigen-decomposition of $\bm{Y}\bm{Y}^\UT$ be $\bm{Y}\bm{Y}^\UT = \bm{U} \mathrm{diag}(\lambda_i) \bm{U}^\UT$, then the action of $\tilde{F}_t$ is defined as $\bm{U} \mathrm{diag}(\tilde{F}_t(\lambda_i)) \bm{U}^\UT$. The other matrix denoisers are defined similarly. The final estimates at iteration $t$ are formed by entry-wise post-processing functions $\phi_{u,t}, \phi_{v,t}$:
\begin{align}\label{eq:Estimates}
\hat{\bm{u}}_t = \phi_{u,t}(\bm{u}_{\leq t}; \bm{a}), \quad
\hat{\bm{v}}_t = \phi_{v,t}(\bm{v}_{\leq t}; \bm{b}).
\end{align}
\end{definition}
Crucially, the matrix denoisers $F_t$ and $G_t$ are required to be \textit{trace-free} in the following sense:
\begin{align}\label{eq:trace free}
\lim_{M\to\infty} \frac{1}{M}\mathrm{Tr}(F_{t}(\bm{Y}\bm{Y}^\UT))&\ac \langle F_{t}(\lambda) \rangle_{\mu} = 0, \\
\lim_{N\to\infty} \frac{1}{N}\mathrm{Tr}(G_{t}(\bm{Y}^\UT\bm{Y})) &\ac \langle G_{t}(\lambda)\rangle_{\tilde{\mu}} = 0,\notag
\end{align}
where the inner products are taken with respect to the limiting spectral measures $\mu$ and $\tilde{\mu}$ respectively. The iterate denoisers are required to be \textit{divergence-free}:
\begin{align}
  &\mathbb{E}[\partial_s f_t(\mathsf{U}_{\leq t-1};\mathsf{A})] = 0, \quad \forall s \in\{ 1,\dots ,t-1\}, \label{eq:div_free_f} \\
 & \mathbb{E}[\partial_s g_t(\mathsf{V}_{\leq t-1};\mathsf{B})] = 0, \quad \forall s \in \{1,\cdots t-1\},\label{eq:div_free_g}
\end{align}
where $(\mathsf{U}_{\leq t-1},\mathsf{V}_{\leq t-1})$ are state evolution random variables to be defined in the following section. Notice that the state evolution random variables and the iterative denoisers are defined in a recursive way, see \cite[Remark 2]{Dudeja2024optimality} for related discussions.

\subsection{State Evolution}

The state evolution (SE) theory of OAMP predicts that the empirical joint distributions of the iterates and signals converge (in the Wasserstein-2 distance) to a sequence of scalar random variables evolving through Gaussian channels, driven by independent $\mathsf{Z}_{u,t}, \mathsf{Z}_{v,t} \sim \mathcal{N}(0,1)$ from true signals:
\begin{align}
(\mathsf{U}_{*}, \mathsf{A}) &\sim \pi_U, \quad \mathsf{U}_t = \mu_{u,t} \mathsf{U}_{*} + \sigma_{u,t}\mathsf{Z}_{u,t}, \label{eq:SEuchannel} \\ 
(\mathsf{V}_{*}, \mathsf{B}) &\sim \pi_V, \quad \mathsf{V}_t = \mu_{v,t} \mathsf{V}_{*} + \sigma_{v,t}\mathsf{Z}_{v,t},
\label{eq:SEvchannel}
\end{align}
where the iterate denoiser outputs are modeled as $\mathsf{F}_t = f_t(\mathsf{U}_{\leq t-1};\mathsf{A})$ and $\mathsf{G}_t = g_t(\mathsf{V}_{\leq t-1};\mathsf{B})$. The parameters of these channels are updated recursively for each iteration $t \in \mathbb{N}$:

\begin{align}
    \mu_{u,t} &= \alpha_t \langle F_t(\lambda) \rangle_{\nu_1} + \beta_t (1+\delta^{-1}) \langle \sigma\widetilde{F}_t(\sigma^2) \rangle_{\nu_3}, \label{eq:mu_u_se} \\
    \mu_{v,t} &= \beta_t \langle G_t(\lambda) \rangle_{\nu_2} + \alpha_t (1+\delta) \langle \sigma\widetilde{G}_t(\sigma^2) \rangle_{\nu_3}, \nonumber \\
    \sigma_{u,t}^2 &= \begin{aligned}[t]
                                & \alpha_t^2 \langle F_t^2(\lambda) \rangle_{\nu_1} + \beta_t^2 \delta^{-1} \langle \lambda \widetilde{F}_t^2(\lambda) \rangle_{\nu_2} + \sigma_{f,t}^2 \langle F_t^2(\lambda) \rangle_{\mu} \\
                                & \quad + 2\alpha_t\beta_t (1+\delta^{-1}) \langle \sigma F_t(\sigma^2)\widetilde{F}_t(\sigma^2) \rangle_{\nu_3} - \mu_{u,t}^2 \\
                                & \quad + \sigma_{g,t}^2 \delta^{-1} \langle \lambda \widetilde{F}_t^2(\lambda) \rangle_{\tilde{\mu}}, \nonumber
                          \end{aligned}  \\
    \sigma_{v,t}^2 &= \begin{aligned}[t]
                                & \beta_t^2 \langle G_t^2(\lambda) \rangle_{\nu_2} + \alpha_t^2 \delta \langle \lambda \widetilde{G}_t^2(\lambda) \rangle_{\nu_1} + \sigma_{f,t}^2 \delta \langle \lambda \widetilde{G}_t^2(\lambda) \rangle_{\mu} \\
                                & \quad + 2\alpha_t\beta_t (1+\delta) \langle \sigma G_t(\sigma^2)\widetilde{G}_t(\sigma^2) \rangle_{\nu_3} - \mu_{v,t}^2 \\
                                & \quad + \sigma_{g,t}^2 \langle G_t^2(\lambda) \rangle_{\tilde{\mu}},\nonumber
                          \end{aligned} 
\end{align}
where the signal alignments $\alpha_t, \beta_t$ and the residual variances of the denoiser outputs, $\sigma^2_{f,t}, \sigma^2_{g,t}$ are defined as:
\begin{align}
\alpha_t &\triangleq \mathbb{E}[\mathsf{U}_*\mathsf{F}_t], & \beta_t &\triangleq \mathbb{E}[\mathsf{V}_*\mathsf{G}_t], \label{eq:se_alignments}\\
\sigma^2_{f,t} &\triangleq \mathbb{E}[\mathsf{F}_t^2] - \alpha_t^2, & \sigma^2_{g,t} &\triangleq \mathbb{E}[\mathsf{G}_t^2] - \beta_t^2. \label{eq:se_variances}
\end{align}
In the above expressions, the inner products are taken with respect to the limiting spectral measures $\nu_1, \nu_2, \nu_3$, whose characterization is provided in Lemma~\ref{lem:spectral_measures}.

The following theorem characterizes the high-dimensional dynamics of the OAMP algorithm; its proof is provided in the full version of this paper \cite[Appendices~B and C]{Chen2025rect}.

\begin{theorem}[State Evolution]\label{Thm:StateEvolution}
Under Assumption~\ref{assump:main}, the iterates of the OAMP algorithm from Definition~\ref{def:OAMP} (initialized with iterates independent of $\bm{W}$) converge in $\mathcal{W}_2$ distance for any finite $t \in \mathbb{N}$ to the state evolution variables of \eqref{eq:SEuchannel}-\eqref{eq:SEvchannel}:
\begin{align}
(\bm{u}_*, \bm{u}_1, \dots, \bm{u}_t; \bm{a}) & \wc (\mathsf{U}_*, \mathsf{U}_1, \dots, \mathsf{U}_t; \mathsf{A}), \\
(\bm{v}_*, \bm{v}_1, \dots, \bm{v}_t; \bm{b}) & \wc (\mathsf{V}_*, \mathsf{V}_1, \dots, \mathsf{V}_t; \mathsf{B}),
\end{align}
where the state evolution random variables $(\mathsf{U}_t,\mathsf{V}_t)_{t\ge1}$ are described in \eqref{eq:SEuchannel}-\eqref{eq:se_variances}.
\end{theorem}

\subsection{Optimal OAMP Algorithm}\label{sec:optimal OAMP}
It is possible to derive the optimal OAMP algorithm directly from the state evolution. The detailed derivation can be found in \cite[Appendix~D]{Chen2025rect}; here, we present only the final form. 

The resulting optimal OAMP algorithm is given by
\begin{align*}
\bm{u_t}^* &= \frac{1}{\sqrt{w_{1,t}}}\big[ \bm{F}^*_{t} \bar{\phi}(\bm{u}_{t-1}^*|w_{1,t-1}) + \bm{\widetilde{F}}^*_{t} \bm{Y}\bar{\phi}(\bm{v}^*_{t-1}|w_{2,t-1})\big], \\
\bm{v_t}^* &= \frac{1}{\sqrt{w_{2,t}}}\big[ \bm{G}^*_{t} \bar{\phi}(\bm{v}_{t-1}^*|w_{2,t-1}) + \bm{\widetilde{G}}^*_{t} \bm{Y}^\UT\bar{\phi}(\bm{u}^*_{t-1}|w_{1,t-1})\big].
\end{align*}
The update matrices are given by $\bm{F}^*_{t} \triangleq F^*_t(\bm{Y}\bm{Y}^\UT)$ and $\bm{\widetilde{F}}^*_{t} \triangleq \tilde{F}^*_t(\bm{Y}\bm{Y}^\UT)$, with analogous definitions for $\bm{G}^*_{t}$ and $\bm{\widetilde{G}}^*_{t}$ applied to $\bm{Y}^\UT\bm{Y}$. The final estimates are given by:
\begin{align}\label{eq:post-process}
\hat{\bm{u}}_t^* \triangleq \phi(\bm{u_t}^*|w_{1,t}),\quad \hat{\bm{v}}_t^* \triangleq \phi(\bm{v_t}^*|w_{2,t}).
\end{align}
The components of the algorithm are defined as follows:
\begin{enumerate}
\item \textit{Iterate Denoisers:} For a scalar Gaussian channel $\mathsf{X} = \sqrt{\omega}\mathsf{X}_* + \sqrt{1-\omega}\mathsf{Z}$, the MMSE estimator for the signal $\mathsf{X}_*$ with side information $\mathsf{C}$ is  $\phi(x;c|\omega) \triangleq \mathbb{E}[\mathsf{X}_* | \mathsf{X}=x, \mathsf{C}=c]$. As required by \textit{divergence-free} condition \eqref{eq:div_free_f}, the DMMSE estimator $\bar{\phi}$ is given by \cite{Dudeja2024optimality}:

\begin{equation}\label{eq:dmmse_estimator_def}
\bar{\phi} (x; c| \omega) \triangleq  
\frac{ \phi(x;c | \omega) - \frac{\mathbb{E}[\mathsf{Z}\phi(\mathsf{X};\mathsf{C} | \omega)]}{\sqrt{1-\omega}} x}{ 1 - \frac{\sqrt{\omega}}{\sqrt{1-\omega}} \mathbb{E}[\mathsf{Z}\phi(\mathsf{X};\mathsf{C} | \omega)]},\quad \forall \omega\in(0,1).
\end{equation}
 \item \textit{Matrix Denoisers:} The \textit{trace-free} denoisers are:
 \begin{subequations}
\begin{align}
F_t^*(\lambda)
&\triangleq
(1+\rho_{1,t}^{-1})\cdot
\left[
1-\langle P_t^*(\lambda)\rangle_{\mu}^{-1}\cdot P_t^*(\lambda)
\right],\notag
\\
\widetilde{F}_t^*(\lambda)
&\triangleq
(1+\rho_{2,t}^{-1})\cdot
\left[
\langle P_t^*(\lambda)\rangle_{\mu}^{-1}\cdot \widetilde{P}_t^*(\lambda)
\right], \label{eq:tFstar}
\\
G_t^*(\lambda)
&\triangleq
(1+\rho_{2,t}^{-1})\cdot
\left[
1-\langle Q_t^*(\lambda)\rangle_{\widetilde{\mu}}^{-1}\cdot Q_t^*(\lambda)
\right],\notag
\\
\widetilde{G}_t^*(\lambda)
&\triangleq
(1+\rho_{1,t}^{-1})\cdot
\left[
\langle Q_t^*(\lambda)\rangle_{\widetilde{\mu}}^{-1}\cdot \widetilde{Q}_t^*(\lambda)
\right].
\label{eq:tGstar}
\end{align}
where
\begin{align}
    P_t^*(\lambda;\rho_{1,t},\rho_{2,t}) &\triangleq \frac{\lambda(\rho_{2,t}\varphi_2(\lambda) + \delta)}{D_t(\lambda)},\label{eq: Pstar} \\
    \widetilde{P}_t^*(\lambda;\rho_{1,t},\rho_{2,t}) &\triangleq \frac{\sqrt{\delta}\rho_{2,t}\varphi_3(\lambda)}{D_t(\lambda)}, \\
    Q_t^*(\lambda;\rho_{1,t},\rho_{2,t}) &\triangleq 
\frac{\delta\lambda(\rho_{1,t}\varphi_1(\lambda) + 1)}{D_t(\lambda)} ,
 \label{eq:Qstar} \\
    \widetilde{Q}_t^*(\lambda;\rho_{1,t},\rho_{2,t}) &\triangleq \frac{\sqrt{\delta}\rho_{1,t}\varphi_3(\lambda)}{D_t(\lambda)},
\end{align}
and 
\begin{equation*}
    D_t(\lambda) \triangleq (\rho_{1,t}\varphi_1(\lambda) + 1)(\rho_{2,t}\varphi_2(\lambda) + \delta)\lambda - \rho_{1,t}\rho_{2,t}\varphi_3^2(\lambda).
\end{equation*}
\end{subequations}
The shrinkage functions $\varphi_1,\varphi_2,\varphi_3$ are defined in \eqref{Eqn:varphi_functions}.
\end{enumerate}

\subsection{State Evolution of Optimal OAMP}
The MSE performance of the optimal OAMP algorithm admits the following simplified form.
\begin{proposition}\label{Pro:OAMP_FP}
For the optimal OAMP algorithm, $\|\hat{\bm{u}}_t-\bm{u}_\ast\|^2/M\xrightarrow[]{\mathrm{a.s.}}\mathrm{mmse}_{\mathsf{U}}(w_{1,t})$ and $\|\hat{\bm{v}}_t-\bm{v}_\ast\|^2/N\xrightarrow[]{\mathrm{a.s.}}\mathrm{mmse}_{\mathsf{V}}(w_{2,t})$, where
\begin{align}
\rho_{1,t} &= \frac{1}{\mathrm{mmse}_{\mathsf{U}}(w_{1,t-1})} - \frac{1}{1-w_{1,t-1}}, \nonumber \\
\rho_{2,t} &= \frac{1}{\mathrm{mmse}_{\mathsf{V}}(w_{2,t-1})} - \frac{1}{1-w_{2,t-1}}, \nonumber\\
 w_{1,t} &= 1 - \frac{1-\langle P^*_t(\lambda) \rangle_{\mu}}{\langle P^*_t(\lambda) \rangle_{\mu} }\cdot \frac{1}{\rho_{1,t}}, \nonumber \\
 w_{2,t} &= 1 - \frac{1-\langle Q^*_t(\lambda)\rangle_{\tilde{\mu}}}{\langle Q^*_t(\lambda) \rangle_{\tilde{\mu}}} \cdot \frac{1}{\rho_{2,t}}, \label{eq:OptimalRecursion}
\end{align}
with $\mathrm{mmse}_\mathsf{X}(w) \triangleq \mathbb{E}\left[(\mathsf{X}_* - \mathbb{E}[\mathsf{X}_*|\mathsf{X}])^2\right]$ for the channel $\mathsf{X} = \sqrt{w}\mathsf{X}_* + \sqrt{1-w}\mathsf{Z}$.
Here, we used the shorthand $P^*_t \triangleq P^*(\cdot;\rho_{1,t},\rho_{2,t})$ in \eqref{eq: Pstar} and similarly for $Q^*_t$ in \eqref{eq:Qstar}.
\end{proposition}

We now specialize our results to the case where the noise matrix has i.i.d.\ Gaussian entries. In this setting, the fixed-point equations of the optimal OAMP algorithm coincide with those of the standard AMP algorithm~\cite{montanari2021estimation}, up to a re-parameterization.

\begin{proposition}\label{thm:OAMP_Wigner_FP}
Suppose $\bm{W}$ has i.i.d. $\mathcal{N}(0, 1/N)$ entries. Then, the fixed-point equations of \eqref{eq:OptimalRecursion} can be simplified to:
\begin{align*}
\mathrm{mmse}_\mathrm{U}(w_1) &= 1 - \frac{1}{\theta^2} \frac{w_2}{1-w_2}, \\
\mathrm{mmse}_\mathrm{V}(w_2) &= 1 - \frac{\delta}{\theta^2} \frac{w_1}{1-w_1}.
\end{align*}
\end{proposition}

The proofs of Propositions \ref{Pro:OAMP_FP} and \ref{thm:OAMP_Wigner_FP} are provided in \cite{Chen2025rect}.

\section{Simulation Results}

\subsection{I.I.D. Gaussian Noise}
\begin{figure}[htb]
\begin{minipage}[b]{.48\linewidth}
  \centering
  \centerline{\includegraphics[width=\linewidth]{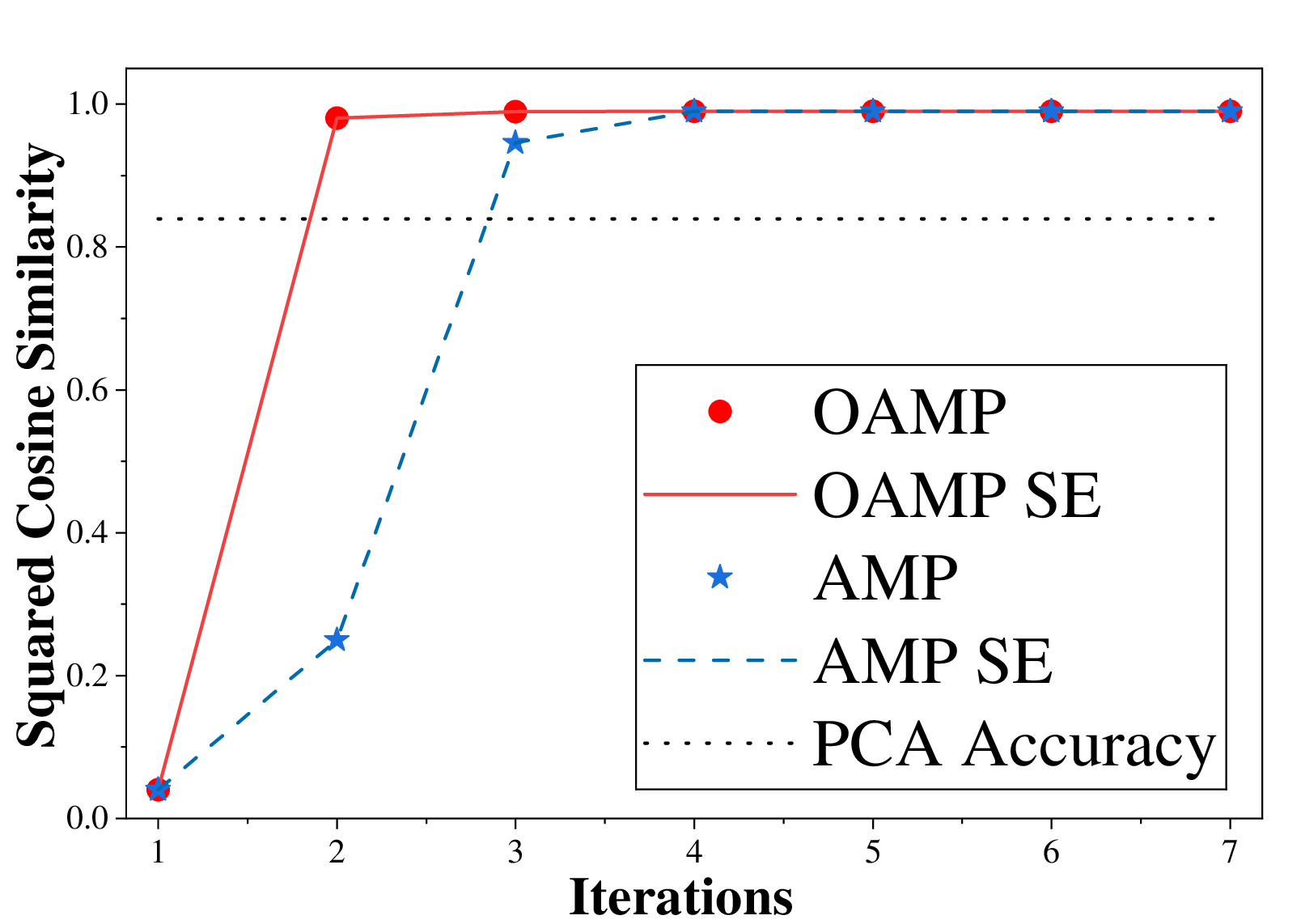}}
  \centerline{(a) $\bm{u}$-channel}\medskip
\end{minipage}
\hfill
\begin{minipage}[b]{0.48\linewidth}
  \centering
  \centerline{\includegraphics[width=\linewidth]{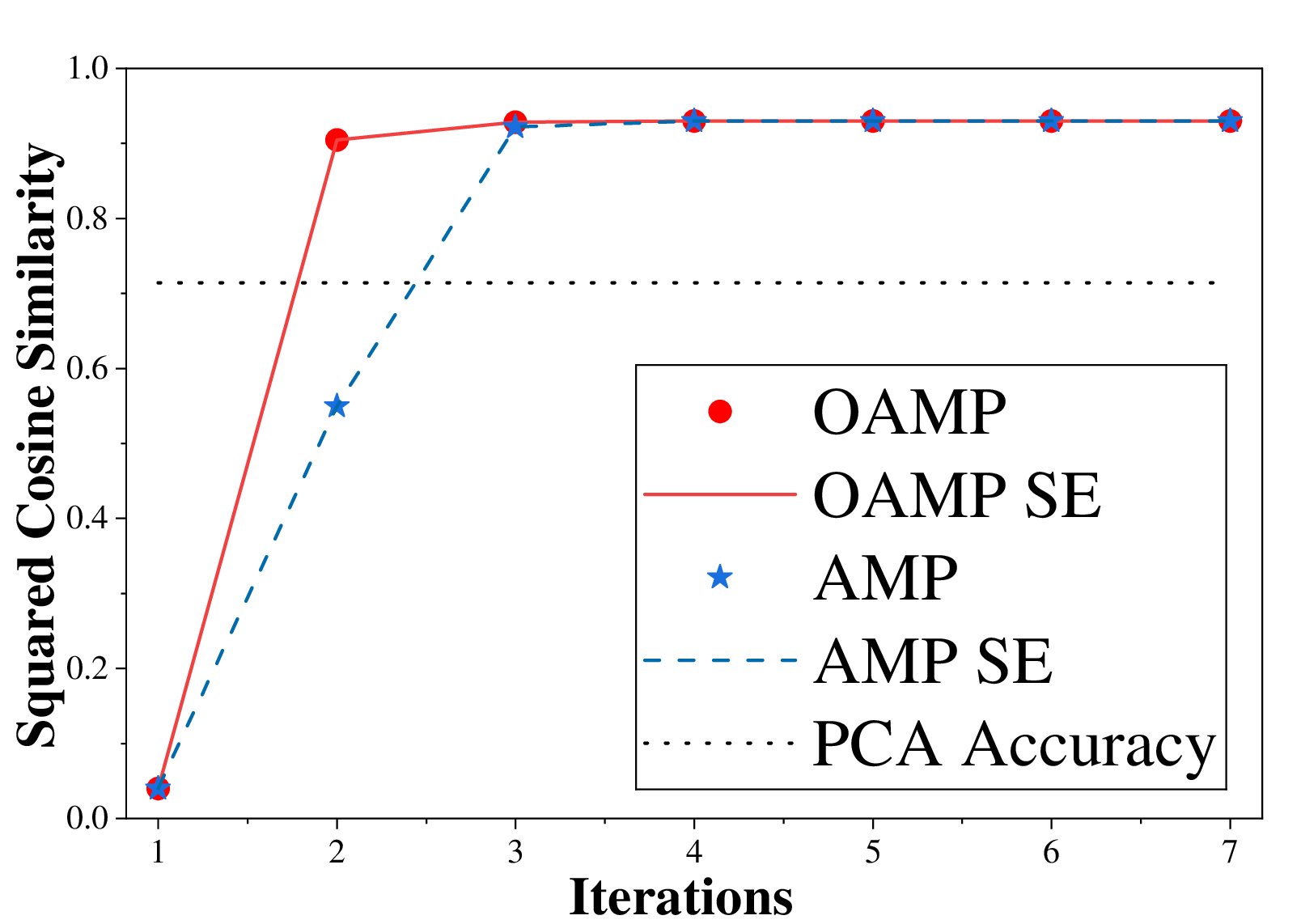}}
  \centerline{(b) $\bm{v}$-channel}\medskip
\end{minipage}
\caption{Simulation results for the i.i.d. Gaussian noise model. The empirical results (markers) are averaged over $50$ runs with $M=4000, N=8000$ and $\theta=2$. The algorithm is initialized with a prior cosine similarity of $0.2$ for both signals.}
\label{fig:res}
\end{figure}

We begin with simulations under i.i.d. Gaussian noise. The priors for the true signals $\bm{u}_*$ and $\bm{v}_*$ are i.i.d. Rademacher. Fig.~\ref{fig:res} reports the squared cosine similarity performances of PCA, AMP, and the proposed OAMP algorithm. The results show that the empirical performance of OAMP agrees closely with its state evolution predictions. Moreover, the proposed OAMP algorithm converges to the same fixed point as AMP, consistent with Proposition~\ref{thm:OAMP_Wigner_FP}.

\begin{figure}[htb]
\begin{minipage}[b]{.48\linewidth}
  \centering
  \centerline{\includegraphics[width=\linewidth]{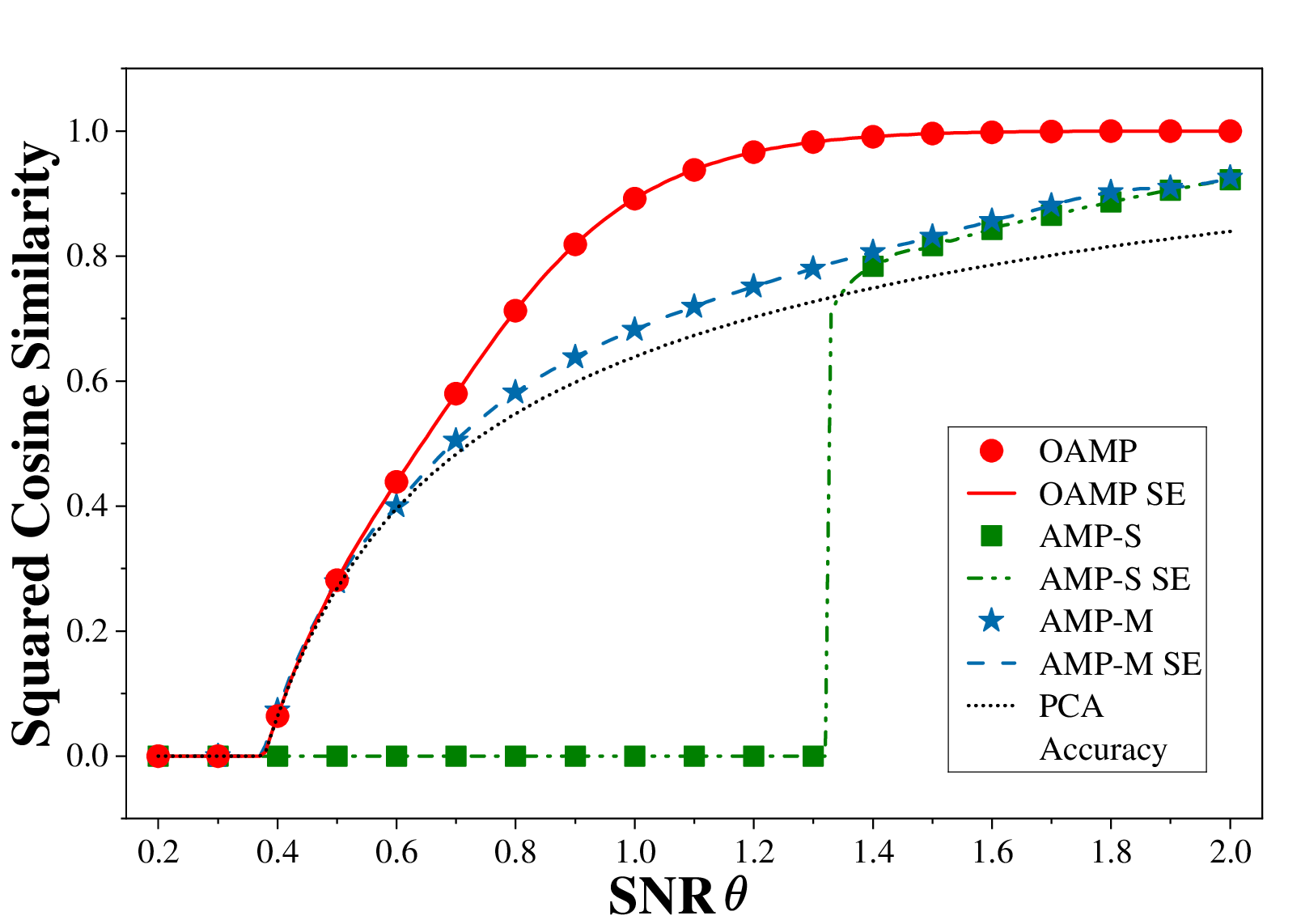}}
  \centerline{(a) $\bm{u}$-channel}\medskip
\end{minipage}
\hfill
\begin{minipage}[b]{0.48\linewidth}
  \centering
  \centerline{\includegraphics[width=\linewidth]{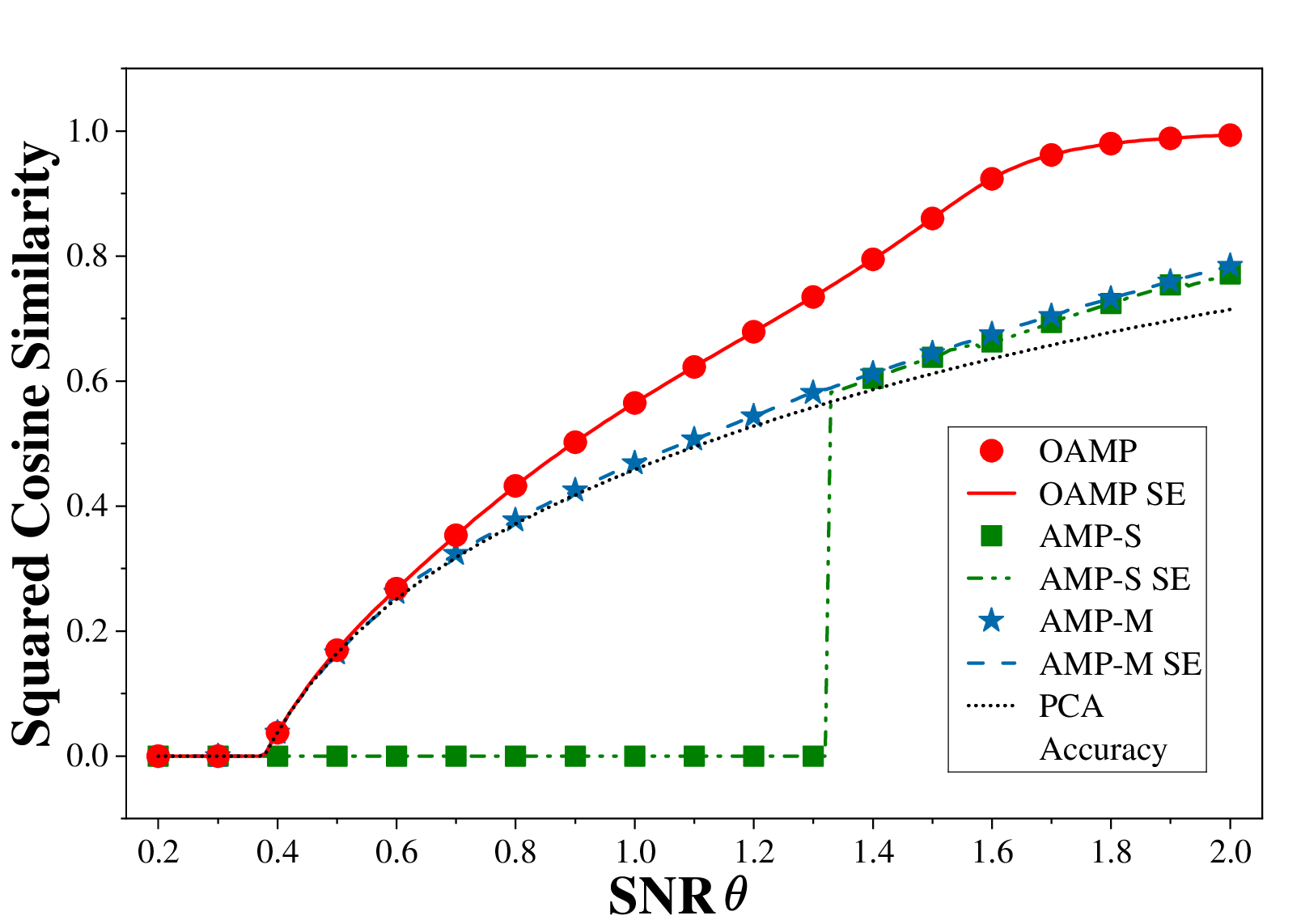}}
  \centerline{(b) $\bm{v}$-channel}\medskip
\end{minipage}
\caption{Simulation results for the non-Gaussian noise case. The eigenvalues of $\bm{WW}^\top$ are drawn from a rescaled/shifted Beta(1.5, 1.5) distribution with support on $[1, 3]$. The signal components are from a Rademacher prior. Other settings are the same as in Fig.~\ref{fig:res}.}
\label{fig:beta_spectrum}
\end{figure}

\subsection{Non-Gaussian Noise}

For the rotationally-invariant noise model, we compare the proposed OAMP with RI-AMP \cite{fan2022ROT}. We evaluate two variants of RI-AMP: a single-iterate version (AMP-S) and a multi-iterate version (AMP-M). In AMP-S, the denoiser is the MMSE function applied to the current iterate, while in AMP-M the denoiser optimally combines all past signal-plus-noise observations using their covariance before applying a MMSE denoiser (see \cite[Remark~3.3]{fan2022ROT} for details). Figure~\ref{fig:beta_spectrum} shows that the optimal OAMP algorithm consistently outperforms PCA and both RI-AMP variants.

\subsection*{Acknowledgments}

This work is supported by the National Key R\&D Program of China under Grant No. 2024YFA1014200, and the National Natural Science Foundation of China under Grant under No. 62571526.

\vfill
\pagebreak

\bibliographystyle{IEEEbib}
\bibliography{strings,refs}

\end{document}